%

\magnification=\magstep1
\input amstex
\documentstyle{amsppt}
\NoBlackBoxes

\def\c{\cite}
\topmatter

\title
Large ordinals
\endtitle

\author
Thomas Jech
\endauthor

\affil
The Pennsylvania State University \endaffil \abstract{Let $j$ be an
elementary embedding of $V_{\lambda}$ into $V_{\lambda}$ that is not
the identity, and let $\kappa$ be the critical point of $j$.  Let $\Cal A$
be the closure of $\{j\}$ under the operation $a (b)$ of application,
and let $\Omega$ be the closure of $\{\kappa\}$ under the 
operation \linebreak
$\min \{\xi: a(\xi) \ge b(\alpha)\}$.

We give a complete description of the set $\Omega$ under an assumption
(Threshold Hypothesis) on cyclic left distributive algebras.}
\endabstract

\endtopmatter

\document

\baselineskip 20pt

\head{\bf 1.  Introduction} \endhead

Let $\lambda$ be a limit ordinal such that there exists a nontrivial 
elementary 
embedding $j: (V_{\lambda}, \epsilon) \to (V_{\lambda}, \epsilon)$.  
The existence 
of such a $\lambda$ is a large cardinal axiom, and by Kunen \c{5}, $\lambda = 
\lim\limits_{n\to\infty} \kappa_n$ where $\kappa_0=\kappa$ is the critical 
point of $j$ and 
$\kappa_{n+1} = j(\kappa_n)$ for all $n=0,1,2,\dots$\,.

Following Laver \c{6}, if $j$ and $k$ are elementary embeddings from 
$V_{\lambda}$ to $V_{\lambda}$, let $j\cdot k$ denote the elementary embedding 
$j (k)= \bigcup_ {\alpha < \lambda}j  (k \upharpoonright V_{\alpha})$.  The 
binary operation $j\cdot k$ satisfies the left-distributive law 
$$
	a (bc) = ab (ac)  \tag{\text{LD}}
$$
[Here and throughout the paper we adopt the convention that $abc=(ab)c$.]

Let $j = V_{\lambda} \to V_{\lambda}$ and let $\kappa = \text{ crit} (j)$
be the critical point of $j$.  Let $\Cal A = \Cal A_j$ be the closure
of $\{j\}$ under $\cdot$ and let $\Cal P = \Cal P_j$ be the closure of
$\{j\}$ under $\cdot$ and $\circ$ where $\circ$ denotes composition.
Let $$\Gamma = \{a\kappa : a \in \Cal A\}.$$  
[Another convention we adopt is
writing $a\kappa$ instead of $a(\kappa)$.]

As $\text{ crit} (ab) = a \text{ crit} (b)$, $\Gamma$ is the set of
all critical points of all $a \in \Cal A$, and it is easily seen that
$\Gamma = \{a\kappa: a \in \Cal P\} = \{\text{crit} (a) : a \in \Cal
P\}$.

In \c{6}, Laver proved that $\Cal A$ is the free left distributive 
algebra on one 
generator, and in \c{7} he showed, using a result of Steel, 
that $\Gamma$ has order type $\omega$.  In fact, if we let 
$$
	j_1 = j, \quad j_{n+1} = j_n j
$$
then $\Gamma = \{\gamma_n: n = 0,1,2,\dots\}$ where 
$$
	\gamma_n = \text{ crit} (j_{2^n}) \quad \text{ and } \quad 
	\gamma_0 < \gamma_1 <\dots< \gamma_n < \dots \;.
$$ 

In this paper we investigate certain ordinal numbers defined in terms 
of the embeddings in $\Cal A $ ($\Cal P)$ and the critical points in
$\Gamma$.  These ordinals have been studied in \c{7}, \c{1} and \c{3}.

\definition
{\bf 1.1.  Definition} \ The set $\Sigma$ of {\it simple ordinals} is
the closure of $\Gamma$ under the operation 
$$
	a"\alpha = \sup \{a \xi : \xi < \alpha\} \qquad\qquad (a \in \Cal A)
$$

The set $\Omega$ of {\it ordinals} is the closure of $\Gamma$ under the operation
$$
	a^-\alpha = \min \{\xi : a \xi \ge \alpha\} \qquad\qquad (a \in \Cal A)
$$
\enddefinition

The following facts are consequences of elementarity:

\proclaim
{\bf 1.2.  Lemma} \ For all $a, b \in \Cal P$ and all ordinals $\alpha$,
$$
\aligned
	a(b\alpha) & = ab (a\alpha) \\
	a(b" \alpha) & = ab" a\alpha \\
	a" b"  \alpha & = (a \circ b)" \alpha \\
	a (b^-\alpha) & = ab^- a\alpha \\
	a^-b^-\alpha & = (b\circ a)^- \alpha.
\endaligned
$$
As a corollary, $\Sigma = \{a" \gamma : \gamma \in \Gamma \text{ and }
a \in \Cal P\}$ and $\Omega = \{a^- \gamma: \gamma \in \Gamma \text{ and } 
a \in \Cal P\}$.
\endproclaim

The following argument shows that every ordinal in $\Sigma$ is in $\Omega$.

\proclaim
{\bf 1.3.  Lemma} (Dougherty) \ Let $c \in \Cal A$ be such that $\gamma = 
\text{ crit } (c)$.  Then $a"\gamma = c^- (c a \gamma)$, for every $a 
\in \Cal P$.
\endproclaim

\demo
{\bf Proof} \ We have $c(a" \gamma) = ca" c\gamma$, and because
$c\gamma > \gamma$, it follows that $ca" c\gamma > ca\gamma$.  Thus 
$c(a"\gamma) > 
ca\gamma$.

If $\eta < a"\gamma$, then $\eta < a\xi$ for some $\xi < \gamma$.
Since $c\xi = \xi$, we have $c\eta < c (a \xi) = ca (c\xi) = ca\xi <
ca \gamma$.  Thus $a"\gamma$ is the least $\eta$ such that $c\eta > c
a \gamma$.

\medpagebreak

It has been conjectured by Laver that $\Sigma=\Omega$.  We prove this
equality (Theorem 3.9), under an assumption on cyclic LD algebras (the
Threshold Hypothesis 3.1).  Under the same hypothesis, we give a
complete description of ordinals in $\Omega$ (Theorem 4.4).

To conclude this introduction, we state the following facts about the
ordinals that will be used in subsequent arguments:
\enddemo

\proclaim
{\bf 1.4.  Lemma} \ (a) \ If $\alpha < \beta$ then $a \alpha < a"\beta$ 
\newline

\hskip .8in \;(b) \ If $\gamma \in \Gamma$ and $a \gamma > \gamma$ then 
$a"\gamma < a \gamma$. \newline

\hskip .8in \;(c) \ If $\gamma \in \Gamma$ and $b\gamma > \gamma$ then 
$a"b \gamma > a (b"\gamma)$
\endproclaim

\demo
{\bf Proof} \ (a) follows from the definition; (b) from the fact that
cf$(a"\gamma) = \gamma$ while cf$(a\gamma) = a (\text{cf}\gamma) = a
\gamma$, and (c) combines (a) and (b).
\enddemo

\head{\bf 2.  Critical points and cyclic LD algebras} \endhead

We shall exploit the remarkable connection between the algebras $\Cal
A$ and $\Cal P$ and the (finite) cyclic left-distributive algebras.
We shall first review some facts from \newline
\c{3} about cyclic LD algebras.

For each $n$ let $A_n = \{0,1,\dots,2^n -1\}$.  There is a unique
left-distributive operation $*_n$ on $A_n$ such that $a*1= a+1 \text{
mod } 2^n$.  For every $a \in A_n$ there exists a number $p_n (a) =
2^k$, the period of $a$ such that
$$
	a < a*1 < a *2 <\cdots < a* (2^k-1), \quad a*2^k = 0   \tag2.1
$$
and $a* (2^k + b) =  a*b$.

In particular,
$$
	p_n (0) = 2^n, \;p_n (2^n -1) = 1, \;p_n (2^{n-1}) = 2^{n-1}
\tag2.2
$$
and for all $a$, if $0< a<2^n-1$ then $1< p_n(a) < 2^n$.

Reduction modulo $2^n$ is a homomorphism from $A_{n+1}$ to $A_n$:
$$      
	a \;*_{n+1} b \text{ mod } 2^n = (a \text{ mod } 2^n) \;*_n 
	(b \text{ mod } 2^n).
\tag2.3
$$ It follows that for every $a \in A_n$, $p_{n+1} (a)$ either remains
equal to $p_n(a)$ or doubles:
$$
	p_{n+1} (a) = \cases \;\;p_n (a) \\  2p_n (a)  \endcases
\tag2.4
$$
and 
$$
	p_{n+1} (a + 2^n) = p_n (a).  \tag2.5
$$

If the period $p_n(a) = 2^k$ doubles to $2^{k+1}$, 
then $a \;*_{n+1} 2^k = 2^n$.

\definition
{\bf 2.1. Definition} \ The {\it threshold} $t_n(a)$ of $a \in A_n$ is 
the least 
$c$ such that \newline  $a \;*_n c \ge \nomathbreak 2^{n-1}$:
$$
	a \;*_n (t_n (a) -1) < 2^{n-1} \le a \;*_n t_n(a).  \tag2.6
$$
We have $t_n (a) \le 2^{k-1}$ where $2^k = p_n (a)$.

By (2.3), the inverse limit of the $A_n$ is a left-distributive algebra; let 
$A_{\infty}$ denote its subalgebra generated by the element 1.

If $w$ is an element of the free left-distributive algebra on one
generator 1 (a ``word"), let $[w]_n$ denote the element of $A_n$ to
which $w$ evaluates.  By (2.3) we have, for all $w$, 
$$
	[w]_{n+1} = \cases [w]_n \\ [w]_n + 2^n \endcases
\tag2.7
$$ 
and $A_{\infty} \vDash v = w$ iff for all $n$, $[v]_n = [w]_n$.  Theorem
4.4 of \c{3} gives several conditions equivalent to the statement that
$A_{\infty}$ is the free algebra; by \c{6}, these are true under the
assumption of the existence of a nontrivial elementary embedding $j :
V_{\lambda} \to V_{\lambda}$.
\enddefinition

If a nontrivial $j: V_{\lambda} \to V_{\lambda}$ exists, then $\Cal A$
is the free one-generated left-distributive algebra and $\Cal A$ is
isomorphic to $A_{\infty}$.  Moreover, as Laver has shown in \c{6},
the equivalence relation $k \overset{\gamma_n}\to = \ell$ (defined in
\c{6}) gives a homomorphism of $\Cal A$ onto $A_n$.  We recall that
this equivalence relation can be defined algebraically on
$A_{\infty}$; see \c{3}, Definition 5.1. In particular we have
$$
	a  \overset{\gamma_n}\to =  b \quad \text{ iff } \quad 
	[a]_n = [b]_n.  \tag2.8
$$
As $\Cal A$ and $A_{\infty}$ are isomorphic, we shall identify the
generator $j$ of $\Cal A$ with 1, and consider the elements $a \in
A_{\infty}$ to be elementary embeddings.  In particular, every integer
$k$ can be identified with some $k \in \Cal A$ via
$$
	1 = j, \qquad k+1 = k * j.   \tag2.9
$$

Using (2.8), we note that if $A_n \vDash a = b$ then for every 
$\gamma \in \Gamma$
$$
\alignedat5
	\text{if } &\quad a\gamma &< \gamma_n &\quad \text{ then } &\quad 
	b\gamma &= a\gamma \\
	\text{if } &\quad a\gamma &\ge \gamma_n &\quad \text{ then } &\quad 
	b\gamma &\ge \gamma_n \\
	\text{if } &\quad a"\gamma &< \gamma_n &\quad \text{ then } &\quad 
	b"\gamma &= 
a" \gamma 
\endalignedat
\tag2.10
$$

\definition
{\bf 2.1. Definition} \ For every word $a$ let $s(a)$ (the {\it signature} 
of $a$) 
be the largest $n$ such that $[a]_n=0$; if $a$ is an integer then $s(a)$ 
is the 
largest $s(a)$ such that $2^{s(a)}$ divides $a$.  

The following summarizes the connection between the algebras $A_n$ and the 
critical points $\gamma \in \Gamma$:
\enddefinition

\proclaim
{\bf 2.3. Lemma} \ $\text{ crit } (a) =  \gamma_{s(a)}$; \newline

\hskip .9in $a \gamma_k \ge \gamma_n \text{ iff } p_n(a) \le 2^k$; \newline

\hskip .9in $a \gamma_k = \gamma_n \text{ where } n=s (a * 2^k)$.

\noindent In particular, this includes Laver's 
result mentioned in the introduction:
$$
	\text{ crit } (2^n)= \gamma_n,\,\, 2^n\gamma_n = \gamma_{n+1}; 
	\tag2.11
$$
the latter is equivalent to this fact about the $A_n$'s:
$$
	A_{n+2} \vDash 2^n * 2^n = 2^{n+1},  \tag2.12
$$
which has the following consequence (that one can also prove directly):
$$
	A_{n+2} \vDash 2^n * a = 2^n +a \qquad\qquad (a \le 2^n).  \tag2.13
$$
\endproclaim

\proclaim
{\bf 2.4. Lemma} \ $(2^n -1) \gamma_0 = \gamma_n$,  $(2^n-1) \gamma_1 = 
\gamma_{n+1}$. 
\endproclaim

\demo
{\bf Proof} \ 
$(2^n-1) \gamma_0 = (2^n-1) \text{ crit } (1) = \text{ crit } ((2^n-1) *1) = 
\text{ crit }  (2^n) = \gamma_n$; $(2^n-1) \gamma_1 = (2^n-1) (1 \gamma_0) = 
(2^n-1) *1 ((2^n-1) \gamma_0) = 2^n \gamma_n = \gamma_{n+1}$
\enddemo

\proclaim
{\bf 2.5.  Lemma} \ If $a < 2^n-1$ then $a \gamma_0 < \gamma_n$.
\endproclaim

\demo
{\bf Proof} \ $s (a*1) < n$.
\enddemo

Following \c{3}, Section 3, let
$$
	a \,\circ_n b = (a \;*_n (b+1)) -1 \text{ mod } 2^n,
$$ 
for all $a, b \in A_n$; the relation $\circ_n$ on $A_n$ is a homomorphic
image of composition on $\Cal P$ under the homomorphism given by the
equivalence relation $=^{\gamma_n}$.

\proclaim
{\bf 2.6. Lemma} \ Let $a \in \Cal P$.  If $a \gamma_0 = \gamma_n$ then 
$a"\gamma_1 = (2^n-1)" \gamma_1$.  If $k\ge 1$ and $a \gamma_k = \gamma_n$ 
then $a"(2^k-1)"\gamma_1 = (2^n-1)"\gamma_1$.
\endproclaim

\demo
{\bf Proof} \ We prove the second statement only, as the proof of the first
one is similar.  We have $a"(2^k-1)"\gamma_1 = a \circ (2^k-1)"\gamma_1$, 
$a \,\circ_{n+1} (2^k-1) = (a \;*_{n+1} 2^k)-1 = 2^n-1$, and the statement 
follows 
by Lemma 2.4 and by (2.10).
\enddemo

\proclaim
{\bf 2.7. Corollary} \  $\gamma_1$ is the least ordinal in $\Omega$ above 
$\gamma_0$; for every $k \ge 1$, $(2^k-1)" \gamma_1$ is the least ordinal in 
$\Omega$ above $\gamma_k$.
\endproclaim

\demo
{\bf Proof} \ Again, we only prove the second statement.  Let $\alpha = 
a^-\gamma$ be an ordinal in $\Omega$ greater than $\gamma_k$.  Let 
$a \gamma_k = \gamma_n$;  since $\gamma_k < \alpha$, we have 
$\gamma_n < \gamma$.  By Lemmas 2.6 and 2.4, $a" (2^k-1)" \gamma_1 = 
(2^n-1)"\gamma_1 < (2^n-1) \gamma_1 = \gamma_{n+1} \le \gamma$ and so 
$(2^k-1)" \gamma_1 \le \alpha$.
\enddemo

\proclaim
{\bf 2.8. Lemma} \ If $2^n < a < 2^{n+1}$ then there is no $\alpha$ such 
that $\gamma_n \le a \alpha < \gamma_{n+1}$, and for no $\alpha$, 
$\gamma_n \le a" \alpha < \gamma_n$.
\endproclaim

\demo
{\bf Proof} \ Let $b = a-2^n$; we have $A_{n+1} \vDash a = 2^n * b$.  If 
$\alpha \ge \gamma_n$ then $a \alpha \ge \gamma_{n+1}$ because 
$\text{ crit } (a) \le \gamma_n$.  If $\alpha < \gamma_n$ then $2^n \alpha 
= \alpha$, and
$$
\align
	a\alpha & = 2^n b (2^n\alpha) = 2^n (b\alpha), \\
	a" \alpha & = 2^n b" 2^n\alpha = 2^n (b" \alpha).
\endalign
$$
Hence both $a\alpha$ and $a"\alpha$ are in the range of $2^n$, which is 
disjoint from the interval $[\gamma_n, \gamma_{n+1})$.
\enddemo

\proclaim
{\bf 2.9. Corollary} \ Every $\alpha \in \Sigma$ between $\gamma_n$ and 
$\gamma_{n+1}$ is equal to $a"\gamma$ for some $\gamma \in \Gamma$ and 
$a < 2^n$.
\endproclaim

\demo
{\bf Proof} \ Let $\alpha = b" \gamma$ where $b \in \Cal P$ and let 
$a = [b]_{n+1}$.  Hence $\alpha = a" \gamma$, and $a < 2^n$ by Lemma 2.8.
\enddemo

\head{\bf 3.  The Threshold Hypothesis and its consequences} \endhead

We shall now formulate a conjecture about the cyclic algebras and use
it to prove results about embeddings in $\Cal P$ and ordinals in
$\Omega$.

\subhead{3.1.  The Threshold Hypothesis (TH)} \endsubhead  Let $a < 2^n-1$, 
and let $p_n (a) = 2^k$.  If $c + 1$ is the threshold of $a$ in $A_n$ then 
$p_{k+1} (c) = 2p_k (c)$.

[We recall that the conclusion of TH is equivalent to the statement that 
$\gamma_k$ is in the range of $c$.]

We conjecture that (TH) holds in every $A_n$.  In the applications that 
follow 
we only use the following consequence of TH:

\proclaim
{\bf 3.2. Lemma} \ Assume TH and let $a < 2^{n-1}$.  If $a \gamma_k = 
\gamma_n$ then there exists a $c < 2^{k-1}$ such that $\gamma_k \in 
\text{ range } (c)$ and $A_n \vDash a c < 2^{n-1}$.
\endproclaim

\demo
{\bf Proof} \ By Lemma 2.3, $p_n (a) = 2^k$, and because $a < 2^{n-1},
k\ge 2$.  Let $c+1$ be the threshold of $a$ in $A_n$; by (2.6), 
$0< c < 2^{k-1}$, and by TH, $\gamma_k \in \text{ range } (c)$.  
Since $c < t_n (a)$, we have 
$A_n \vDash a c < 2^{n-1}$.  By Lemma 2.3, $a \gamma_k = \gamma_n$
implies that $t_{n+1} (a) = 2^k$, and it follows that $A_{n+1} \vDash a c < 
2^{n-1}$.
\enddemo

\proclaim
{\bf 3.3. Theorem} \ Assume TH.  If $a < 2^n$ and $a \gamma_{k+1} = 
\gamma_{n+1}$ then $a" \gamma_{k+1} > \gamma_n$.
\endproclaim

\demo
{\bf Proof} \ The statement is vacuously true for $n=0$.  For every $n\ge 1$ 
we prove the theorem by downward induction on $a$.

First let $a = 2^n-1$.  By Lemma 2.4 we have $k=0$, $a\gamma_0 = 
\gamma_n$ and $a \gamma_1 = \gamma_{n+1}$.  Therefore $\gamma_n < a" 
\gamma_1 < \gamma_{n+1}$.

Now let $a < 2^n-1$.  by Lemma 3.2 there exists some $c<2^k$ such that 
$A_{n+1} \vDash a c < 2^n$, and $\gamma_{k+1} \in \text{ range } (c)$.  Let 
$\delta \in \Gamma$ be such that $c \delta = \gamma_{k+1}$, and let 
$b = [ac]_{n+1}$.  Since $c < 2^k$, we have (by Lemmas 2.3 and 2.5) 
$\gamma_0 < \delta < \gamma_{k+1}$ and so $a \delta = \gamma_{m+1}$ for some 
$m < n$. Therefore
$$
	ac \gamma_{m+1} = ac (a \delta) = a (c \delta) = a \gamma_{k+1} = 
	\gamma_{n+1}.
$$
In $A_{n+2},$  $ ac = b$ or $ac = b+2^{n+1}$.  The latter is impossible
because $A_{n+2} \vDash b +2^{n+1} = 2^{n+1} b$ (by (2.13)),
$2^{n+1}
\gamma_{m+1} = \gamma_{m+1}$ (because $m< n$), and so $2^{n+1} b 
\gamma_{m+1} = 2^{n+1} b (2^{n+1} \gamma_{m+1}) = 2^{n+1} (b \gamma_{m+1}) 
\ne \gamma_{n+1}$, by (2.11).  Hence $[ac]_{n+2} = b$ and so $b \gamma_{m+1} = 
\gamma_{n+1}$.  Since $b>c$ and $b<2^n$, we have, by the induction hypothesis, 
$b" \gamma_{m+1} > \gamma_n$.  Now the statement for $a$ follows (using 
Lemma 1.4 and (2.10)):
$$
	a" \gamma_{k+1} = a" c \delta > a (c" \delta) = ac" a\delta = 
	ac" \gamma_{m+1} = b" \gamma_{m+1} > \gamma_n.
$$
\enddemo

\proclaim
{\bf 3.4. Corollary} (TH) \ If $a < 2^n$ and if $\gamma_n < a" \gamma < 
\gamma_{n+1}$ then $a\gamma = \gamma_{n+1}$
\endproclaim

\demo
{\bf Proof} \ Let $a \gamma = \gamma_{m+1}$.  By Theorem 3.3, $a"\gamma > 
\gamma_m$, and so $m=n$.
\enddemo

\proclaim
{\bf 3.5. Corollary} (TH) \ If $a < 2^n$ and $\gamma_{n+1} \in 
\text{ range } (a)$ then there exists an $\alpha \in \Sigma$ such that 
$\gamma_n < a \alpha < \gamma_{n+1}$.
\endproclaim

\demo
{\bf Proof} \ It is $\alpha = c" \delta$ in the proof of Theorem 3.3.
\enddemo

\definition
{\bf 3.6. Definition} \ For $\alpha \in \Sigma$, let $\alpha^+$ denote
the least ordinal in $\Sigma$ greater than $\alpha$.
\enddefinition

\proclaim
{\bf 3.7. Corollary} (TH) \ Let $a < 2^n$ and let $\alpha^+ = 
\gamma \in \Gamma$ be a critical point.  If $\gamma_n \le a \alpha < 
\gamma_{n+1}$ then $a\gamma = \gamma_{n+1}$.
\endproclaim

\demo
{\bf Proof} \ If $\gamma = \gamma_1$ then $\alpha = \gamma_0$ by 
Corollary 2.7, and so $a \gamma_0 = \gamma_n$.  By Lemmas 2.4 and 2.5, 
$a = 2^n-1$ and so $a\gamma_1 = \gamma_{n+1}$.

Thus let $\gamma = \gamma_{k+1}$ where $k>0$.  Then $a \gamma_{k+1} = 
\gamma_{m+1}$ where $m\ge n$.  By Corollary 3.5 there exists a $\beta \in 
\Sigma$ such that $\gamma_m < a \beta < \gamma_{m+1}$, and since $\gamma = 
\alpha^+$, we have $\gamma_m < a \alpha$.  It follows that $m=n$ and hence 
$a \gamma = \gamma_{n+1}$.
\enddemo

\proclaim
{\bf 3.8. Corollary} (TH) \ If $\gamma = \alpha^+ \in \Gamma$ and if 
$a \in \Cal P$ is arbitrary, then there is no critical point $\delta \in 
\Gamma$ between $a \alpha$ and $a" \gamma$.

\endproclaim

\demo
{\bf Proof} \ This is true if $a\alpha = \alpha$, so assume that 
$\text{ crit } (a) < \alpha$.
Let $n$ be such that $\gamma_n \le a \alpha < \gamma_{n+1}$.  
We will show that $a" \gamma = \gamma_{n+1}$.  By Lemma 2.8, $A_{n+1} 
\vDash a < 2^n$.  Let $b = [a]_{n+1}$; then $b < 2^n, \gamma_n \le b \alpha 
< \gamma_{n+1}$, and by Corollary 3.7 we have $\gamma = \gamma_{n+1}$.

Since $\text{ crit } (b) = \text{ crit } (a) < \gamma_n$ we have $b"
\gamma < \gamma_{n+1}$, and it follows that $a" \gamma < \nomathbreak
\gamma_{n+1}$.
\enddemo

\proclaim
{\bf 3.9. Theorem} \ Assume TH.  Then Laver's Conjecture holds; i.e. every 
ordinal $u^-\lambda \in \Omega$ is in $\Sigma$.
\endproclaim

\demo
{\bf Proof} \ Let $\xi = u^-\lambda (\lambda \in \Gamma)$ be a
counterexample.  Let $b" \gamma \, (\gamma \in \Gamma)$ be the least
ordinal in $\Sigma$ greater than $\xi$.  Let $\alpha \in \Sigma$ be
such that $\alpha^+ = \gamma$.  Since $b \alpha \in \Sigma$, we have
$b\alpha < \xi < b" \gamma$.

Let $a = u \circ b$.  We have 
$$
	a" \gamma = (u \circ b)" \gamma = u" b"\gamma > u \xi \ge \lambda
$$
(by Lemma 1.4), and 
$$
	a \alpha = (u \circ b) \alpha = u (b\alpha) < \lambda.
$$
This contradicts Corollary 3.8.

\medpagebreak

We conclude this Section with the following observation of Laver:
\enddemo

\proclaim
{{\bf 3.10. Lemma} (Laver)} \ If $\Omega = \Sigma$ then $$a" \alpha^+ = 
(a \alpha)^+$$ for all $a \in \Cal P$ and all $\alpha \in \Sigma$.
\endproclaim

\demo
{\bf Proof} \ Assume otherwise and let $\xi \in \Sigma$ be such that 
$a \alpha < \xi < a" \alpha^+$.  Then $\alpha < a^- \xi < \alpha^+$, a 
contradiction.
\enddemo

\head{\bf 4.  Ordinals between {$\bold \gamma_n$} and {$\bold
\gamma_{n+1}$}} \endhead

In this Section we again assume TH and describe all $\Omega$-ordinals
between two consecutive critical points.  We also formulate the {\it
Uniqueness Hypothesis}, another conjecture about the cyclic algebras,
and use it to prove that the representation is unique.

\definition
{\bf 4.1. Definition} \ Let $\alpha \in \Sigma$ be such that $\gamma_n < 
\alpha < \gamma_{n+1}$; $\alpha$ is {\it special} (below $\gamma_{n+1}$) if 
$\alpha^+ = \gamma_{n+1}$.  
\enddefinition

\proclaim
{\bf 4.2. Lemma} \ If $\alpha$ is special then there exist no 
$\xi \in \Omega$, 
$\xi < \alpha$ and no $a \in \Cal P$ such that $\alpha = a \xi$.
\endproclaim

\demo
{\bf Proof} \ If $a \xi = \alpha$ then by Lemma 3.9 the ordinal $\alpha^+ = 
a"\xi^+$ has cofinality cf~$\xi^+ < \alpha^+$ and therefore is not a critical 
point.
\enddemo

\proclaim
{\bf 4.3. Lemma} \ If $\alpha \in \Omega$, $\gamma_n < \alpha < 
\gamma_{n+1}$ and $\alpha$ is not special, 
then there exist  $\xi \in \Omega$, 
$\xi < \alpha$ and $a \in \Cal P$ such that $\alpha = a \xi$.
\endproclaim

\demo
{\bf Proof} \ If $\alpha$ is not special then  $\alpha^+ = a"\gamma$ 
for some $a \in \Cal P$ and $\gamma \in \Gamma$, $\gamma \ne \gamma_0$.  Let 
$\xi \in \Omega$ be such that $\xi^+ = \gamma$.  By Lemma 3.10, $(a \xi)^+ = 
a" \xi^+ = a" \gamma = \alpha^+$, and so $\alpha = a\xi$.
\enddemo

\proclaim
{\bf 4.4. Theorem} \ Let $a<2^n$, $\alpha \in \Omega$ and $\gamma_n < \alpha < 
\gamma_{n+1}$.  The ordinal $\alpha$ is not special if and only if there exist
$c < 2^n, \gamma \in \Gamma, a < c, \lambda \in \Gamma$ and $b$ such that 
$\alpha = c" \gamma, a\lambda = \gamma$ and $A_{n+1} \vDash ab=c$.
\endproclaim

\demo
{\bf Proof} \ First assume that $\alpha$ is not special.  By Lemma 4.3
there exist $\xi < \alpha$ and $a$ such that $\alpha = a \xi$.  By
Lemma 2.8 we may assume that $a < 2^n$.

Let $b$ and $\lambda \in \Gamma$ be such that $\xi = b" \lambda$.
Then $\alpha = a \xi = a (b" \lambda) = ab" a \lambda$.  Let $a
\lambda = \gamma$ and $c = [ab]_{n+1}$.  We have $a<c,$  $c"\gamma =
ab" a\lambda = \alpha$, and by Lemma 2.8, $c<2^n$.

Conversely, assume that the condition holds.  Then $\alpha = c"\gamma
= ab" a\lambda = a(b"\lambda)$.  As $a<2^n$, its critical point is
below $\gamma_n$, hence $a\alpha > \alpha$ and so $b"\lambda < 
\alpha$.  By Lemma 4.2, $\alpha$ is not special.
\enddemo

We not describe all $\Omega$-ordinals between consecutive critical
points:

\proclaim
{\bf 4.5. Theorem} \ Let $n>0$.  There exist a finite sequence
$$
	a_1 = 2^n-1 > a_2 >\cdots> a_{k_n}
$$
of embeddings and a finite sequece 
$$
	\gamma_{i_1} = \gamma_1,\dots,\gamma_{i_{k_n}}
$$
of critical points such that
$$
	\gamma_n < a"_1 \gamma_{i_1} < \cdots < a"_{k_n} \gamma_{i_{k_n}} < 
	\gamma_{n+1}
$$
are all the $\Omega$-ordinals between $\gamma_n$ and $\gamma_{n+1}$, and 
for every $k=1,\dots,k_n-1$,
$$
	a_k" \gamma_{i_k} = a_{k+1} \alpha
$$
where $\alpha$ is the special ordinal below $\gamma_{i_{k+1}}$.
\endproclaim    

\demo
{\bf Proof} \ Let $a_1 = 2^n-1$.  By Corollary 2.7, $a"_1 \gamma_1$
is the least $\Omega$-ordinal greater than $\gamma_n$ (this is proved
without using TH).  By induction, let $k \ge 1$, $a_k = a$ and
$\gamma_{i_k} = \gamma$.  If $a"\gamma$ is special below
$\gamma_{n+1}$ we are done.  Thus assume that $a" \gamma$ is not
special, and so there exist a $\xi < a"\gamma$ and some $x$ such that 
$x\xi = a"\gamma$.

Let $\alpha$ be the least $\xi$ such that $x\xi = a"\gamma$ for some
$x$.  We claim that $\alpha$ is special: if not then $\alpha =
y"\eta$ for some $y$ and $\eta < \alpha$; then $(x\circ y) \eta =
x(y\eta) = x\alpha = a"\gamma$, contradicting the minimality of
$\alpha$.  Therefore there exists a special $\alpha < a"\gamma$ such
that $c\alpha = a"\gamma$.  By Corollary 2.9, there is such a $c$
with the property that $c<2^n$.

Let $a_{k+1}$ be the largest $c<2^n$ such that for some special
$\alpha, a"\gamma \le c\alpha < \gamma_{n+1}$.  First we note that
by the induction hypothesis, it is impossible that $c\ge a$: this is
clear for $k=1$, and if $k>1$, then this would contradict the fact
that $a_k$ is the largest $a < 2^n$ such that for some special $\xi$,
$a"_{k-1} \gamma_{i_{k-1}} \le a \xi < \gamma_{n+1}$.  Hence $c<a_k$.

We conclude the proof by showing that $c\alpha = a"\gamma$.  Thus 
assume that $c \alpha > a"\gamma$.  There exist a $b$ and some 
$\lambda \in \Gamma$ such that $a=b"\gamma$, and we have 
$$
	c\alpha = c(b" \lambda) = cb" c\lambda = d" c\lambda
$$
where $d = [cb]_{n+1}$.  By Lemma 2.8 we have $c < d < 2^n$.

Let $\eta$ be special below $c\lambda$.  By Lemma 3.10, $d" c\lambda$
is the successor of $d\eta$, and so $d\eta \ge a"\gamma$.  This
contradicts the maximality of $c$.
\enddemo 

We shall now address the question of {\it uniqueness} of the
representation given by Theorem 4.5.  If $a"\gamma = b"\delta$ and
$\gamma, \delta \in \Gamma$ then, by the reason of cofinality, $\gamma
= \delta$.  The question is whether we can have $a" \gamma = b"
\gamma$ when $a, b<2^{n-1}$ and $a" \gamma < \gamma_{n+1}$.  We prove 
the uniqueness of $a"\gamma$ under the assumption of the following 
{\it Uniqueness Hypothesis}:

\subhead{4.6.  The Uniqueness Hypothesis (UH)} \endsubhead  \ Let 
$a, b < 2^{n-1}$ and let $p_n(a) = p_n(b) = 2^k$.  Let $c$ be the 
least $c$ such that $\gamma_k$ is in the range of $c$ and let 
$c\gamma_i = \gamma_k$.  If $a\gamma_i = b\gamma_i$ and $A_n 
\vDash ac = bc$, then $a=b$.

We conjecture that (UH) holds in every $A_n$.  The proof of Theorem 4.8 
uses the following consequence of UH:

\proclaim
{4.7. Lemma} \ Assume UH and let $a, b <2^{n-1}$, $a \ne b$.  If $a
\gamma_k = b\gamma_k = \gamma_n$ and if $c$ is the least $c$ such that
$\gamma_k \in \text{ range }(c)$, then $A_n \vDash ac\ne bc$.
\endproclaim

\demo
{\bf Proof} \ Let $a, b < 2^{n-1}$ and let $c$ be least with
$\gamma_k \in \text{ range }(c)$; let $c \gamma_i = \gamma_k$.  Assume
that $A_n \vDash ac = bc$.

By TH, $c$ is smaller than the threshold of either $a$ or $b$ in $A_n$
and so $[ac]_n = [bc]_n < 2^{n-1}$.  Since $\gamma_n$ is in the range
of both $a$ and $b$, $[ac]_{n+1} = [bc]_{n+1} < 2^{n-1}$.  Thus
$\gamma_n = a (c\gamma_i) = ac (a\gamma_i) = [ac]_{n+1} (a\gamma_i) =
[bc]_{n+1} (b\gamma_i)$, and it follows that $a\gamma_i = b\gamma_i$.
By UH, $a=b$.
\enddemo

\proclaim
{\bf 4.8. Theorem} \ Assume UH.  If $a, b <2^n$ and $\gamma_n <
a"\gamma = b"\gamma < \gamma_{n+1}$ then $a=b$.
\endproclaim

\demo
{\bf Proof} \ We proceed by induction on $n$ and, for a given $n$, by
downward induction on $a<2^n$.  If $a = 2^n-1$ then $\gamma = 
\gamma_1$, and because $\gamma_1 = \gamma^+_0$, Lemma 3.10 implies that 
if $b"\gamma_1 = a"\gamma_1$ then $b\gamma_0 = a\gamma_0 =
\gamma_n$, and by Lemma 2.5 we have $b = 2^n-1$.

Now let $b < a < 2^n-1$ be such that $\gamma_n < a"\gamma = b"\gamma
< \gamma_{n+1}$.  By Corollary 3.4, $a\gamma = b\gamma =
\gamma_{n+1}$.  Let $c" \delta$ be the special ordinal below
$\gamma$.  By Theorem 4.5 and the induction hypothesis on $n$, $c$ is
the least $c$ such that $\gamma \in \text{ range }(c)$, and so by
Lemma 4.7, $A_{n+1} \vDash ac \ne bc$.  By TH and by Lemma 2.3 we have 
$[ac]_{n+1} < 2^n$ and $[bc]_{n+1} < 2^n$.  By Lemma 3.10, $a"\gamma 
= (a(c"\delta))^+$, and it follows that $[ac]_{n+1}" a\delta = 
[bc]_{n+1}" b\delta$.  This contradicts the induction hypothesis on 
$a$, since $a < [ac]_{n+1} < 2^n$.
\enddemo 

\head{\bf 5.  The conjectures TH and UH} \endhead

The main result of our paper depends on the conjecture TH for finite
LD algebras.  In this Section we discuss the numerical evidence for
the conjecture as well as related conjectures.

The statement TH is formulated in terms of the algebras $A_n$.  In
principle, one can verify the validity of the statement TH for any
particular value of $n$.  In practice, the number of calculations for
$A_n$ grows exponentially, so we can't really expect to verify TH for
too large values of $n$.

In our experiments we use a sophisticated software developed by
Randall Dougherty.  Using various unpublished results about the
algebras $A_n$, Dougherty devised an algorithm that can compute the
binary operation in the algebras $A_n$ for all $n \le 48$.  We used
the code \c{2} with the author's permission.  I am grateful to the
Computer Science Department of Penn State for letting me use their
equipment that is vastly superior to the computers available to
mathematicians at Penn State who don't happen to specialize in
differential equations.

\subhead{5.1.  Experimental result} \endsubhead  TH is true for all 
$n \le 30$.

The proof of Lemma 3.2 that uses TH uses only the instance of TH when 
$\gamma_n$ is in the range of the embedding $a$.  Thus we can 
formulate a weaker hypothesis that is still sufficient for the results 
of our paper.

\subhead{5.2. The weak threshold hypothesis (WTH)} \endsubhead  
Let $a < 2^{n-1}$ be such that $a\gamma_k = \gamma_n$.  If $c+1$ is
threshold of $a$ in $A_n$ then $\gamma_k$ is in the range of $c$.

WTH can be reformulated in several equivalent ways.  To see that, we 
first observe the following:

\proclaim
{5.3.  Lemma} \  Let $a < 2^{n-1}$ be such that $\gamma_n \in 
\text{ range }(a)$ and let $c+1 = t_n (a)$.  Then
$$
[ac]_{n+1}  = [ac]_n = [ac]_{n-1}, $$
$$[a\circ c]_{n+1} = [a\circ c]_n.
$$
\endproclaim

\demo
{\bf Proof} \ As $c+1 = t_n (a)$, we have $[ac]_n < 2^{n-1} \le 
[a(c+1)]_n$.  Thus $[ac]_n = [ac]_{n-1}$.  Since $\gamma_n \in  
\text{ range } (a)$, we have $t_{n+1}(a) = p_n(a) > c+1$, and so 
$[ac]_{n+1} = [ac]_n$,  and 
$[a\circ c]_{n+1} = [a(c+1)-1]_{n+1} = [a(c+1)-1]_n = [a\circ c]_n$.
\enddemo  

\proclaim
{5.4.  Lemma} \  Let $a < 2^{n-1},$  $a\gamma_k = \gamma_n$ and 
$c+1 = t_n (a)$.  Then the following are equivalent.

(i) \ $\gamma_k \in \text{ range } (c)$ 

(ii) \ $\gamma_n  \in \text{ range } ([ac]_n)$ 

(iii) \ $\gamma_n \in \text{ range } ([ac]_{n-1})$ 

(iv) \ $\gamma_n \in \text{ range } ([a\circ c]_n)$ 

(v) \ $\gamma_{n-1} \in \text{ range } ([a\circ c]_{n-1})$ 
\endproclaim

\demo
{\bf Proof} \ To show that (i) is equivalent to (ii) (and to (iii)), 
first let $c \delta = \gamma_k$; we have
$$
	\gamma_n = a\gamma_k = a(c\delta) = ac(a\delta) = [ac]_{n+1}
	(a\delta) = [ac]_n (a\delta).
$$
Conversely, if $\gamma_n \in \text{ range } ([ac]_n) = 
\text{ range } ([ac]_{n+1})$ then $a\gamma_k = 
\gamma_n \in \text{ range } (ac)$ and so $\gamma_k \in 
\text{ range } (c)$.

To show the equivalence of (i) and (iv), let first $c\delta = \gamma_k$; 
we have 
$$
	\gamma_n = a\gamma_k = a(c\delta) = (a\circ c) \delta = 
	[a\circ c]_{n+1} \delta = [a\circ c]_n \delta.
$$
Conversely, if $a\gamma_k = \gamma_n = [a\circ c]_n \delta = 
[a\circ c]_{n+1} \delta = (a\circ c) \delta = a(c\delta)$, then 
$\gamma_k = c\delta$.

Finally, to show the equivalence of (iv) and (v), we first observe that 
if \newline 
$[a(c+1)]_n = 2^{n-1}$, we have $[a\circ c]_{n-1} = [a\circ c]_n = 
2^{n-1}-1$ and both $\gamma_{n-1}$ and $\gamma_n$ are in the range of 
$2^{n-1}-1$.

If $[a(c+1)]_n > 2^{n-1}$, then
$$
	[a\circ c)_n =  [a\circ c]_{n-1} +2^{n-1} = 2^{n-1} \;
	*_n [a\circ c]_{n-1}
$$
and because $2^{n-1} \gamma_{n-1} = \gamma_n$, we have $\gamma_{n-1} 
\in \text{ range } ([a\circ c]_{n-1})$ if and only if \newline 
$\gamma_n 
\in \text{ range } ([a\circ c]_n)$.
\enddemo

It follows from Lemma 5.4 that (WTH) is equivalent to any of the 
following three statements:

\subhead{5.5 } \endsubhead

(WTH1) \ If $a<2^n, \gamma_n \in \text{ range } (a)$ and $c+1 = t_n(a)$ 
then $\gamma_n\in \text{ range } ([a c]_n)$.

(WTH2) \ If $a<2^n, \gamma_n \in \text{ range } (a)$ and $c+1 = t_n(a)$ 
then $\gamma_n\in \text{ range } ([a\circ c]_n)$.

(WTH3) \ If $a<2^n, \gamma_n \in \text{ range } (a)$ and $c+1 = t_n(a)$ 
then $\gamma_{n-1}\in \text{ range } ([a\circ \nomathbreak c]_{n-1})$.

We remark that the assumption that $\gamma_n \in \text{ range } (a)$ is 
necessary in each (WTH1), (WTH2) and (WTH3):

\demo
{\bf 5.6. Example} \ (i) Let $n=5$, $a=5$; then $c=1$, $ac=6$ and 
$\gamma_5$ is not in the range of $6$.

(ii)  Let $n=10$, $a=34$; then $c=4$, $[a\circ c]_9 = 242$ and 
$\gamma_9$ is not in the range of $242$.
\enddemo

When investigating the weak threshold hypothesis, we notice that in
most cases it is true for trivial reasons, namely because if
$\gamma_n$ is in the range of $a$ then  $\gamma_{n-1}$ is in the range of $a$
as well.  This leads to the following conjecture.  

\subhead{5.7.  The Twin Hypothesis} \endsubhead 
If $n$ is odd, $a<2^{n-1}$ and $\gamma_n \in \text{ range } (a)$, then 
$\gamma_{n-1} \in \text{ range } (a)$.

If $a$ satisfies the Twin Hypothesis then it satisfies WTH:  this is 
because $t_n(a) = 2^{k-1}$ where $a\gamma_k = \gamma_n$, 
and $\gamma_k \in \text{ range } (2^{k-1}-1)$.

\subhead{5.8.  Experimental result} \endsubhead 
The Twin Hypothesis is true for all odd $n \le 31$.

If the Twin Hypothesis holds for an odd $n$ then it holds for $n+1$, 
for all $a$ such that  $2^{n-1} \le a < 2^n$.  We use the following:

\proclaim
{5.9.  Lemma} \  If $a < 2^{2m}$ then $\gamma_{2m+1} 
\in \text{ range }(a)$ if and only if $\gamma_{2m+2} \in 
\text{ range }(2^{2m}+a)$.
\endproclaim

\demo
{\bf Proof} \ We have  $2^{2m} \gamma_{2m} = \gamma_{2m+1}$, and by a
result of Dougherty (unpublished) and Dr\'apal \c{4}, 
$2^{2m} \gamma_{2m+1} = \gamma_{2m+2}$. 
We also have $A_{2m+3} \vDash 2^{2m} a = 2^{2m} + a$.  If $a\gamma_k = 
\gamma_{2m+1}$ then $\gamma_{2m+2} = 2^{2m} \gamma_{2m+1} = 2^{2m} 
(a\gamma_k) = 2^{2m} a (2^{2m} \gamma_k) = 2^{2m} a \gamma_k = 
(2^{2m} + a) \gamma_k$.

Conversely, if $\gamma_{2m+2} \in \text{ range }(2^{2m} + a)$ then 
$2^{2m} \gamma_{2m+1} \in \text{ range } (2^{2m}a)$ and so 
$\gamma_{2m+1}\in \text{ range }(a)$.
\enddemo

\proclaim
{\bf 5.10.  Corollary}  \ If $n$ is odd and the Twin Conjecture holds 
then for every $a$ such that $2^{n-1} \le a < 2^n$, if $\gamma_{n+1} 
\in \text{ range }(a)$ then $\gamma_n \in  \text{ range } (a)$.
\endproclaim

\demo
{\bf Proof} \ For $a = 2^{n-1}$, we have $\gamma_n \in \text{ range
}(a)$ because $2^{n-1} \gamma_{n-1} = \gamma_n$.  Thus let $2^{n-1} <
a < 2^n$, and let $b = a - 2^{n-1}$.  If $\gamma_{n+1} \in \text{
range } (b+2^{n-1})$ then $\gamma_n \in \text{ range } (b)$ by Lemma
5.9, and $\gamma_{n-1} \in \text{ range } (b)$ by the Twin Hypothesis.
Since $\gamma_n = 2^{n-1} \gamma_{n-1}$ we have $\gamma_n \in \text{
range } (2^{n-1} b)$ and so $\gamma_n \in \text{ range } (a)$.
\enddemo

\medpagebreak

Now we turn our attention to the Uniqueness Hypothesis.  When we apply 
(UH) (in Lemma 4.7) we only use a weaker version:

\subhead{5.11} \endsubhead  If $a, b < 2^{n-1}$ are such that 
$a\gamma_k = b\gamma_k = \gamma_n$ and if $c$ is the least $c$ such 
that $\gamma_k \in \text{ range }(c)$, then $[ac]_n = [bc]_n$ implies 
$a=b$.

We have verified both UH and 5.11 for a large number of embeddings:

\subhead{5.12.  Experimental result} \endsubhead  UH is true for all 
$n\le 17$;  5.11 is true for all $n \le 24$.

As a final remark, we observe that the 5.11 does not necessarily hold 
when $c$ is not the least $c$:

\demo
{\bf 5.13.  Example} \ Let $n=9$, $a=48$, $b=192$ and $c=51$.  Then 
$a\gamma_7 = b\gamma_7 = \gamma_9$, $c\gamma_3 = \gamma_7$ and 
$[ac]_9 = [bc]_9 = 243$.
\enddemo

\head{\bf 6. Appendix. Ordinals below $\gamma_{12}$}\endhead

\magnification=\magstep0
\baselineskip 8pt

$$\gamma_{0}$$
$$\gamma_{1}$$
$$1"\gamma_{1}$$
$$\gamma_{2}$$
$$3"\gamma_{1}$$
$$\gamma_{3}$$
$$7"\gamma_{1}$$
$$4"\gamma_{3}$$
$$3"\gamma_{2}$$
$$2"\gamma_{2}$$
$$1"\gamma_{2}$$
$$\gamma_{4}$$
$$15"\gamma_{1}$$
$$12"\gamma_{3}$$
$$3"\gamma_{3}$$
$$\gamma_{5}$$
$$31"\gamma_{1}$$
$$28"\gamma_{3}$$
$$19"\gamma_{3}$$
$$16"\gamma_{5}$$
$$15"\gamma_{2}$$
$$2"\gamma_{3}$$
$$\gamma_{6}$$
$$63"\gamma_{1}$$
$$60"\gamma_{3}$$
$$51"\gamma_{3}$$
$$48"\gamma_{5}$$
$$15"\gamma_{3}$$
$$\gamma_{7}$$
$$127"\gamma_{1}$$
$$124"\gamma_{3}$$
$$115"\gamma_{3}$$
$$112"\gamma_{5}$$
$$79"\gamma_{3}$$
$$64"\gamma_{7}$$
$$63"\gamma_{2}$$
$$50"\gamma_{3}$$
$$48"\gamma_{6}$$
$$47"\gamma_{2}$$
$$34"\gamma_{3}$$
$$32"\gamma_{6}$$
$$31"\gamma_{2}$$
$$18"\gamma_{3}$$
$$16"\gamma_{6}$$
$$15"\gamma_{4}$$
$$14"\gamma_{2}$$
$$1"\gamma_{3}$$
$$\gamma_{8}$$
$$255"\gamma_{1}$$
$$252"\gamma_{3}$$
$$243"\gamma_{3}$$
$$240"\gamma_{5}$$
$$207"\gamma_{3}$$
$$192"\gamma_{7}$$
$$63"\gamma_{3}$$
$$48"\gamma_{7}$$
$$15"\gamma_{5}$$
$$\gamma_{9}$$
$$511"\gamma_{1}$$
$$508"\gamma_{3}$$
$$499"\gamma_{3}$$
$$496"\gamma_{5}$$
$$463"\gamma_{3}$$
$$448"\gamma_{7}$$
$$319"\gamma_{3}$$
$$304"\gamma_{7}$$
$$271"\gamma_{5}$$
$$256"\gamma_{9}$$
$$255"\gamma_{2}$$
$$242"\gamma_{3}$$
$$240"\gamma_{6}$$
$$47"\gamma_{3}$$
$$32"\gamma_{7}$$
$$15"\gamma_{6}$$
$$\gamma_{10}$$
$$1023"\gamma_{1}$$
$$1020"\gamma_{3}$$
$$1011"\gamma_{3}$$
$$1008"\gamma_{5}$$
$$975"\gamma_{3}$$
$$960"\gamma_{7}$$
$$831"\gamma_{3}$$
$$816"\gamma_{7}$$
$$783"\gamma_{5}$$
$$768"\gamma_{9}$$
$$255"\gamma_{3}$$
$$240"\gamma_{7}$$
$$15"\gamma_{7}$$
$$\gamma_{11}$$
$$2047"\gamma_{1}$$
$$2044"\gamma_{3}$$
$$2035"\gamma_{3}$$
$$2032"\gamma_{5}$$
$$1999"\gamma_{3}$$
$$1984"\gamma_{7}$$
$$1855"\gamma_{3}$$
$$1840"\gamma_{7}$$
$$1807"\gamma_{5}$$
$$1792"\gamma_{9}$$
$$1279"\gamma_{3}$$
$$1264"\gamma_{7}$$
$$1039"\gamma_{7}$$
$$1024"\gamma_{11}$$
$$1023"\gamma_{2}$$
$$1010"\gamma_{3}$$
$$1008"\gamma_{6}$$
$$815"\gamma_{3}$$
$$800"\gamma_{7}$$
$$783"\gamma_{6}$$
$$768"\gamma_{10}$$
$$767"\gamma_{2}$$
$$754"\gamma_{3}$$
$$752"\gamma_{6}$$
$$559"\gamma_{3}$$
$$544"\gamma_{7}$$
$$527"\gamma_{6}$$
$$512"\gamma_{10}$$
$$511"\gamma_{2}$$
$$498"\gamma_{3}$$
$$496"\gamma_{6}$$
$$303"\gamma_{3}$$
$$288"\gamma_{7}$$
$$271"\gamma_{6}$$
$$256"\gamma_{10}$$
$$255"\gamma_{4}$$
$$14"\gamma_{3}$$

\magnification=\magstep1

\enddocument
\end